\documentclass{article}
\usepackage{amsfonts}

\begin{document}

\def\det{\mbox{det}}
\def\Spin{\mbox{Spin}}
\def\Tr{\mbox{Tr \-}}
\def\inj{\mbox{inj \-}}
\def\R{{\mathbb R}}
\def\C{{\mathbb C}}
\def\H{{\mathbb H}}
\def\Ca{{\mathbb Ca}}
\def\Z{{\mathbb Z}}
\def\N{{\mathbb N}}
\def\Q{{\mathbb Q}}
\def\Ad{\mbox{Ad \-}}
\def\k{{\bf k}}
\def\l{{\bf l}}
\def\diag{\mbox{diag \-}}
\def\a{{\alpha}}

\author{Ya.~V.~Bazaikin}
\title{On the causal properties of symmetric Lorentzian spaces
\thanks{The author was supported by the Russian Foundation of Basic Research (grants 09-01-00142-a and 09-01-12130-ofi-m), the Leading Scientific Schools grant NSh-7256.2010.1
and the joined project of SB RAS and UrB RAS N 46}}

\date{}

\maketitle

\begin{abstract}

The causal properties of Lorentzian symmetric spaces are investigated in
the paper. There is proved the global hyperbolicity of the Cahen--Wallach Lorentzian
symmetric spaces, corresponding to solvable Lie algebras.

\end{abstract}

\sloppy

\section[]{Introduction}

Symmetric Lorentzian spaces play a significant role in the foundations of
modern supergravity theories \cite{Jose, Baer}. In \cite{Cahen}, a
classification of all Lorentzian symmetric spaces was given. Besides
the universal covers of multi-dimensional de Sitter and anti-de Sitter
spaces, the only simply-connected symmetric spaces are the
Cahen--Wallach spaces $CW_n(A)$, which are constructed in the
following way \cite{Jose}.
Consider the Euclidean space $\R^n$ with coordinates $(x^1, \ldots,
x^n)$ and the plane $\R^2$ with coordinates $(\xi,\eta)$. Let
$A=(a_{ij})$ be a symmetric matrix. Then $CW_n(A)=(\R^2 \times \R^n,
ds_{n,A}^2)$, where
$$
ds_{n,A}^2 = 2d\eta \left( d\xi + \sum_{i,j=1}^n a_{ij} x^i x^j d
\eta \right) +  (dx^1)^2 + \ldots + (dx^n)^2.
$$
One of the main results of the paper is the following

\vskip0.2cm

{\bf Theorem 1.} {\it The space $CW_n(A)$ is globally hyperbolic. }

\vskip0.2cm

{\bf Remark 1.} The spaces $CW_n(A)$ correspond to solvable Lie algebras
in the classification of \cite{Cahen}. It can be shown in the remaining
cases that de Sitter space is globally hyperbolic (for instance,
\cite{Baer}), and anti-de Sitter space is not globally hyperbolic
\cite{O'Neill}.

Theorem 1 turns out to be an immediate consequence of the following
more general result.

\vskip0.2cm

{\bf Theorem 2.} {\it Let $(M, g_1)$ be a Riemannian $C^1$-manifold
and let there exist a point $p_0$ with the property that the exponential
map $\exp_{p_0}: T_{p_0} M \rightarrow M$ is a diffeomorphism. Denote by
$\rho(p)=\rho(p_0,p)$ the Riemannian distance function to the point
$p_0$.

Consider the following Lorentzian metric $g_2$ on $N=\R^2 \times M$:
$$
g_2 = 2 d\eta \left( d\xi -F d\eta \right) + g_1,
$$
where $\xi, \eta$ are the coordinates on $\R^2$ and $F:M \rightarrow \R$
is a continuous function. Assume that
$$
|F(p)| \leq c_1^2 + c_2^2\rho(p)^2
$$
for some constants $c_1, c_2$. Then $(N,g_2)$ is globally
hyperbolic. }

\vskip0.2cm

{\bf Remark 2.} One can relax the conditions of Theorem 2: instead of
the function $\rho$ one can choose the distance function to some subset $S
\subset M$; and the exponential map condition can be substituted by the following one: the restriction of the function $\rho$ to every piecewise
$C^1$-smooth curve in $M$ is a piecewise $C^1$-smooth real function.
The last property is closely related to the cut locus structure of
$M$.

\section[]{Proof of the main results}

Let $(N,g)$ be a time-oriented Lorentzian manifold. The manifold $(N,g)$
is called {\it strongly causal} if for each point $p \in N$ and for
each neighborhood $U$ of $p$ there exists an open neighborhood $V
\subset U$ of $p \in V$ such that each causal curve in $N$, starting
and ending in $V$, is entirely contained in $U$. As usual, a {\it
causal future (past)} $J^+(p)$ ($J^-(p)$) of the point $p$ is the set of
all such points $q$ that either $p=q$ or there exists a future
directed (past directed) causal curve starting at $p$ and ending at
$q$. A strongly causal Lorentzian manifold $(N,g)$ is called {\it
globally hyperbolic} if for each pair $p,q \in N$ the set $J^+(p)
\cap J^-(q)$ is compact.

It is proved in \cite{Bernal1, Bernal2} that global hyperbolicity of
the time-oriented Lorentzian manifold $(N,g)$ is equivalent to the
following condition: $(N,g)$ is isometric to the direct product $\R
\times S$ with the metric $-\beta dt^2 + g_t$, where $\beta$ is a smooth
positive function, $g_t$ is a Riemannian metric on $S$ depending
smoothly on $t$, and each ''section'' $\{ t \} \times S$ is a smooth
space-like Cauchy hypersurface in $M$.

We say that a space-time $(N_1,g_1)$ is {\it asymptotically dominated}
by a space-time $(N_2,g_2)$ if there exists a continuous map $f:N_1
\rightarrow N_2$, satisfying the following two properties:

1) for each causal curve $\gamma$ in $N_1$ the curve $f(\gamma)$ is
causal in $N_2$;

2) for each point $p \in N_2$ the set $f^{-1}(p)$ is compact.

\vskip0.2cm

{\bf Lemma 1}. {\it Let $(N_1,g_1)$ be asymptotically dominated by
$(N_2, g_2)$. Then if $(N_2,g_2)$ is globally hyperbolic, then $(N_1,
g_1)$ is also globally hyperbolic.

}

\vskip0.2cm

{\bf Proof}. Since $(N_2,g_2)$ is globally hyperbolic, there
exists a time function $T': N_2 \rightarrow \R$ \cite{Bernal1,
Bernal2}. Put $T=T' \circ f$. The property 1) of the map $f$ implies that
$T: N_1 \rightarrow \R$
 is a time function on $N_1$. Then $(N_1,g_1)$ is stably causal
and, consequently, it is strongly causal \cite{Beem}.

Let us prove that for each compact set $K' \subset N_2$ the set
$K=f^{-1}(K') \subset N_1$ is compact. Introducing Riemannian
metrics on $N_1$ and $N_2$, we can assume that $N_1$ and $N_2$ are
metric spaces with countable topology bases. So, it suffices
for us to prove that from each sequence $p_i \in K$ we can extract a
convergent subsequence. Therefore, in view of compactness of $K'$, we
can assume that $f(p_i) \rightarrow q \in K'$. This means that there
exits a subsequence $p_i$ converging to the set $K_0=f^{-1}(q)$, which is
compact by the property 2) of the map $f$. Consider an open
$\varepsilon$-neighborhood $U$ of the compact $K_0$ with the compact closure
$K_1=\bar{U}$. Then we can assume that $p_i \in U$ and, passing to a subsequence, we obtain that $p_i \rightarrow p \in K_1$. So far as
$p_i$ converges to $K_0$, we see that $p \in K_0 \subset K$.

Further, if $x \in J^+(p) \cap J^-(q)$, then $f(x) \in J^+(f(p)) \cap
J^-(f(q))$. Therefore, $J^+(p) \cap J^-(q)$ is a subset of the compact
set $f^{-1} (J^+(f(p)) \cap J^-(f(q)))$. So far as the set $J^+(p) \cap
J^-(q)$ is closed in a strongly causal space-time \cite{Beem}, we can
conclude that it is compact. This argument completes the proof of the global
hyperbolicity of $(N_1, g_1)$. The lemma is proved.

\vskip0.2cm

Consider the following metric on $N=\R^3$ with coordinates $(\eta, \xi,
\tau)$:
$$
g = 2 d\eta \left( d\xi - f(\tau) d \eta \right) + d\tau^2,
$$
where $f:\R \rightarrow \R$ is a continuous function. Let
$(N,g)$ be oriented in time by the coordinate vector field $\frac{\partial}{\partial
\eta}$.

\vskip0.2cm

{\bf Theorem 3.} {\it Let there exist constants $c_1, c_2$ such that
$|f(\tau)| \leq c_1^2\tau^2 + c_2^2$. Then the metric $g$ is globally
hyperbolic. }

\vskip0.2cm

{\bf Remark 3.} In the paper \cite{Bazaikin} the global hyperbolicity of
the metric $g$ was proved in the case when $f(\tau)\sim \tau$. However, the
method that was used there does not work in the general case.

{\bf Proof}. Let $f_0(\tau)=c_1^2\tau^2 + c_2^2$. Let us consider the metric
$$
\tilde{g} = 2 d\eta \left( d\xi - f_0(\tau) d \eta \right) +
d\tau^2,
$$
which dominates $g$, and let us show its global hyperbolicity.

First of all, we construct a global time function $T=T(\eta,\xi,\tau)$.
We will seek $T$ in the form:
$$
T=\eta-\Phi(\xi,\tau),\eqno{(1)}
$$
for some smooth function $\Phi:\R^2 \rightarrow \R$. It is easy to
calculate the gradient of $T$ with respect to the metric $\tilde{g}$:
$$
\nabla T = \left( 1 - 2f_0(\tau) \Phi_\xi \right) \partial_\xi
-\Phi_\xi \partial_\eta -\Phi_\tau \partial_\tau.
$$
We will try to find a time function $T$, whose gradient is time-like, that
is
$$
|\nabla T|^2 = 2f_0(\tau)\Phi_\xi^2+\Phi_\tau^2 - 2 \Phi_\xi <0,
$$
or, equivalently,
$$
\Phi_\xi^2 + \frac{\Phi_\tau^2}{2f_0(\tau)} <
\frac{\Phi_\xi}{f_0(\tau)}.\eqno{(2)}
$$
Let $\phi,\psi: \R \rightarrow \R$ be some smooth functions. More
precisely, we will try to find the function $\Phi$ in the following
form:
$$
\Phi(\xi,\tau)=\phi\left(\frac{\xi}{\psi(\tau)}\right).
$$
Then the condition (2) changes into the following inequality:
$$
\phi'\left(\frac{\xi}{\psi(\tau)}\right) \left[ 1+
\frac{\psi'(\tau)^2}{2 f_0(\tau)} \left( \frac{\xi}{\psi(\tau)}
\right)^2\right] < \frac{\psi(\tau)}{f_0(\tau)}.\eqno{(3)}
$$
Take $\varepsilon>0$ and put $\psi(\tau)=\varepsilon
f_0(\tau)=\varepsilon (c_1^2\tau^2 +c_2^2)$, that is
$\psi(\tau)/f_0(\tau)=\varepsilon$. Then
$$
\frac{\psi'(\tau)^2}{2 f_0(\tau)} = \frac{2 \varepsilon^2
c_1^4\tau^2}{c_1^2 \tau^2 + c_2^2} \leq 2 \varepsilon^2 c_1^2.
$$
Consequently, the inequality (3) holds everywhere if we require $\phi$ to satisfy the
following equation:
$$
\phi'(x)\left( 1+ 2 \varepsilon^2 c_1^2 x^2
\right)=\frac{\varepsilon }{2}.
$$
An elementary integration yields
$$
\phi_\varepsilon (x) = \frac{1}{2\sqrt{2}c_1} \arctan \left(
\sqrt{2}\varepsilon c_1 x \right).
$$
So we have proved that the everywhere defined smooth function
$$
T_\varepsilon (\xi,\eta, \tau)=\eta-\frac{1}{2\sqrt{2}c_1} \arctan
\left(
 \frac{\sqrt{2} \varepsilon c_1 \xi}{c_1^2\tau^2 +c_2^2} \right)
$$
is a time function with the time-like gradient. In particular, this
implies strong causality of the metric $g$ \cite{Beem}.

Consider two points $p_i=(\xi_i,\eta_i,\tau_i)$, $i=1,2$ in $N$,
where $T(p_1)<T(p_2)$ (in the opposite case the intersection of the cones $K=J^+(p_1) \cap J^-(p_2)$ is empty). Note that for all $\eta_0,
\xi_0$ the transformation
$$
(\xi,\eta,\tau) \mapsto (\xi+\xi_0,\eta+\eta_0, \tau)
$$
is an isometry of the Lorentzian manifold $(N,g)$. Consequently, we can assume, without
loss of generality,  that $\eta_1=0$.

\vskip0.2cm

{\bf Lemma 2.} {\it Assuming the conditions described above, if $\eta_2 <
\frac{\pi}{4\sqrt{2}c_1}$, then $K=J^+(p_1) \cap J^-(p_2)$ is
compact.

}

\vskip0.2cm

{\bf Proof of Lemma 2.} In the further calculations we put
$$
x=\frac{ \xi}{c_1^2\tau^2 +c_2^2}, x_i=\frac{ \xi_i}{c_1^2\tau_i^2
+c_2^2}, i=1, 2.
$$
Now choose $\varepsilon>0$ small enough so that
$$
\left|\phi_\varepsilon( x_i)\right| \leq \frac{\pi}{4\sqrt{2}c_1} -
\eta_2, i=1,2.
$$
Consider some point $p \in K$ and two causal curves: $\gamma(T)$,
starting at $p_1$ and ending at $p$, and $\delta(T)$, starting at
$p$ and ending at $p_2$ (it is evident that we can choose the time
function $T$ as a regular parameter on these curves). Joining these
curves, we can assume that we have one piecewise causal curve
$\gamma(T)$, starting at $p_1$, going through $p$ and ending at
$p_2$. Then, using monotonicity of the functions $T$ (strict) and $\eta$
(non-strict) along $\gamma$, we have
$$
\phi_\varepsilon (x(T)) \leq \phi_\varepsilon (x(T)) -
\phi_\varepsilon (x_1) + \frac{\pi}{4\sqrt{2}c_1} -\eta_2 < \eta(T)
- \eta_1 + \frac{\pi}{4\sqrt{2}c_1} -\eta_2\leq
$$
$$
\leq \frac{\pi}{4\sqrt{2}c_1}.
$$
The definition of the function $\phi_\varepsilon$ immediately implies that
there exists a constant $d$ such that $x(T)\leq d$ for all $T$. Since the initial point $p$ was taken to be arbitrary, we have $x(p)\leq
d$ for all $p \in K$. Analogously,
$$
\phi_\varepsilon (x(T)) \geq \phi_\varepsilon (x(T)) -
\phi_\varepsilon (x_2) - \frac{\pi}{4\sqrt{2}c_1} +\eta_2 > \eta(T)
- \eta_2 - \frac{\pi}{4\sqrt{2}c_1} +\eta_2\geq
$$
$$
\geq - \frac{\pi}{4\sqrt{2}c_1}.
$$
Therefore, we can assume that for all $p \in K$ the following
inequality takes place:
$$
|x(p)| \leq d.\eqno{(4)}
$$
Further, from the relation (1) we have:
$$
1=\frac{d\eta}{dT} - \phi_\varepsilon'(x) \frac{dx}{dT},
$$
that is
$$
\frac{dx}{dT} \geq - \frac{2}{\varepsilon} \left( 1+ 2 \varepsilon^2
c_1^2 x^2 \right) \geq - \frac{2}{\varepsilon} \left( 1+ 2
\varepsilon^2 c_1^2 d^2 \right).
 $$
Consequently,
$$
\frac{d x}{dT} = \frac{1}{f_0} \frac{d \xi}{dT} -\frac{\xi}{f_0}
\frac{2c_1^2\tau}{f_0} \frac{d\tau}{dT} = \frac{1}{f_0} \frac{d
\xi}{dT} -x \frac{2c_1^2\tau}{f_0} \frac{d\tau}{dT} \geq -
\frac{2}{\varepsilon} \left( 1+ 2 \varepsilon^2 c_1^2 d^2 \right).
$$
This immediately implies
$$
\frac{d \xi}{dT} \geq -f_0 \frac{2}{\varepsilon} \left( 1+ 2
\varepsilon^2 c_1^2 d^2 \right) + 2 x c_1^2\tau \frac{d\tau}{dT}
\geq -f_0 \frac{2}{\varepsilon} \left( 1+ 2 \varepsilon^2 c_1^2 d^2
\right) - 2 d_1 c_1^2 \left| \tau \frac{d\tau}{dT} \right|
$$
Further, we use the causality condition for the curve $\gamma$:
$$
0 \geq 2 \frac{d\eta}{dT} \left( \frac{d\xi}{dT} - f_0
\frac{d\eta}{dT} \right) +\left( \frac{d \tau}{dT} \right)^2 \geq
$$
$$
2 \frac{d\eta}{dT} \left( -f_0 \frac{2}{\varepsilon} \left( 1+ 2
\varepsilon^2 c_1^2 d^2 \right) - 2 d c_1^2 \left| \tau
\frac{d\tau}{dT} \right| - f_0 \frac{d\eta}{dT} \right) +\left(
\frac{d \tau}{dT} \right)^2\eqno{(5)}
$$
Let us assume, first, that $\frac{d\eta}{dT} \geq 1$. Then, continuing
the previous inequality, we obtain
$$
- 2 f_0 \left(1 + \frac{2}{\varepsilon} \left( 1+ 2 \varepsilon^2
c_1^2 d^2 \right) \right) \left( \frac{d\eta}{dT} \right)^2  - 4 d
c_1^2 \frac{d\eta}{dT} \left| \tau \frac{d\tau}{dT} \right|  +\left(
\frac{d \tau}{dT} \right)^2 \leq 0
$$
In view of our assumption about the growth of $\eta$, we can change
the parameter  $T$ on the curve $\gamma$, taking the coordinate $\eta$ as a new
parameter. In this case (dividing by $f_0$), we have:
$$
\frac{1}{f_0} \left( \frac{d \tau}{d\eta} \right)^2 - 4 d c_1^2
\left| \frac{\tau}{f_0} \frac{d\tau}{d \eta} \right| - 2 \left(1 +
\frac{2}{\varepsilon} \left( 1+ 2 \varepsilon^2 c_1^2 d^2 \right)
\right) \leq 0.
$$
So far as
$$
\frac{|\tau|}{\sqrt{f_0}} =
\frac{1}{\sqrt{c_1^2+\frac{c_2^2}{\tau^2}}} \leq \frac{1}{c_1},
$$
we have
$$
\left( \frac{1}{\sqrt{f_0}}\frac{d \tau}{d\eta} \right)^2 - 4 d c_1
\left| \frac{1}{\sqrt{f_0}} \frac{d\tau}{d \eta} \right| - 2 \left(1
+ \frac{2}{\varepsilon} \left( 1+ 2 \varepsilon^2 c_1^2 d^2 \right)
\right) \leq 0
$$
The last estimate implies
$$
\left| \frac{1}{\sqrt{f_0}}\frac{d \tau}{d\eta} - 2d c_1 \right|
\leq \sqrt{4d^2 c_1^2+2 \left(1 + \frac{2}{\varepsilon} \left( 1+ 2
\varepsilon^2 c_1^2 d^2 \right) \right)}.
$$
By virtue of
$$ \frac{1}{\sqrt{f_0}}\frac{d \tau}{d\eta} = \frac{d
}{d\eta} \left( \frac{1}{c_1} \mbox{ arcsinh } \left( \frac{c_1
\tau}{c_2}\right) \right),
$$
we can integrate the last inequality and, using a priori boundedness of
$\eta$, we obtain that the total increment of the function $\tau$ on those
subsegments of $\gamma$, where $\frac{d\eta}{dT} \geq 1$, is bounded
by some constant $D$ depending only on $c_1, c_2, \eta_2$ and on the constants
$\xi_1, \xi_2, d, \varepsilon$ which are determined by the choice of the points $p_1,
p_2$.

If $\frac{d\eta}{dT} \leq 1$, the inequality (5) can be immediately rewritten as:
$$
\left( \frac{d \tau}{dT} \right)^2  - 4 d c_1^2 \left| \tau
\frac{d\tau}{dT} \right|  - 2 f_0 \left(1 + \frac{2}{\varepsilon}
\left( 1+ 2 \varepsilon^2 c_1^2 d^2 \right) \right) \leq 0.
$$
As in the previous case, dividing by $f_0$ and integrating, we obtain that
the total increment of $\tau$ on those subsegments, where
$\frac{d\eta}{dT} \leq 1$, is also universally bounded by the constant
$D$.

Thus we have proved that the coordinate $\tau$ is universally bounded along
each causal curve. The estimate (4) implies that the coordinate $\xi$ is
also universally bounded. Therefore, the set $K$ is contained in a bounded
domain in $\R^3$, that is, the closure of $K$ is compact. Strong
causality implies the closeness and, consequently, compactness of $K$
\cite{Beem}. The lemma is proved.

\vskip0.2cm

Now let $p_1$ and $p_2$ be the same as above but without
any restriction on $\eta_2$. As above, we consider the intersection of the cones
$K=J^+(p_1) \cap J^-(p_2)$. Take a constant $C>0$ and consider
the following Lorentzian metric on the manifold $N'=\R^3(\xi',\eta',\tau')$:
$$
\tilde{g}'^2 = 2 d\eta' \left( d\xi' - f_0'(\tau) d \eta' \right) +
d\tau'^2,
$$
where
$$
f_0'(\tau')= \frac{c_1^2}{C^2} \tau'^2 + \frac{c_2^2}{C^2}=c_1'^2
\tau'^2 +c_2'^2.
$$
Now choose $C$ large enough for the inequality
$$
\eta_2 < \frac{\pi}{4\sqrt{2}c_1'}\eqno{(6)}
$$
to hold. Consider the transformation $\sigma: N \rightarrow N'$:
$$
\sigma(\xi, \eta, \tau) = \left(\frac{\xi}{C^2}, \eta,
\frac{\tau}{C} \right).
$$
It is evident that the transformation of the Lorentzian spaces $\sigma:
(N,ds^2) \rightarrow (N',ds'^2)$ is a conformal diffeomorphism.
Consequently, it maps homeomorphically the intersection $K$ of the causal
future and past of the points $p_1, p_2$ to the intersection $K'$ of
the corresponding causal future and past of the points $p_1'=(\xi_1/C^2, 0,
\tau_1/C)$ and $p_2'=(\xi_2/C^2, \eta_2, \tau_2/C)$. Inequality (6)
and Lemma 2 imply that $K'$ is compact and, therefore, $K$ is also
compact. Our theorem is proved.

\vskip0.2cm

{\bf Remark 4.} The asymptotical growth of the function $f_0$ in Theorem 3
is optimal in the class of power functions.  Indeed, if, in
the conditions of Theorem 3, we consider $f = \tau^{2+\varepsilon}$ for $\varepsilon>0$ , then it is sufficient to take a class of the causal curves
of the form  $\gamma(s) = (\xi_0,\eta(s), \tau(s))$. Such curves are
causal on the plane with coordinates $(\eta, \tau)$ with respect to the
Lorentzian metric
$$
- 2 d \eta^2 + \frac{d\tau^2}{f(\tau)}.
$$
An elementary integration shows that this metric is not globally
hyperbolic: there exist light-like curves, both future- and past-directed, which escape to $+\infty$ with respect to the variable $\tau$ for a finite increment of $\eta$.

\vskip0.2cm

{\bf Proof of Theorem 2.} It suffices to note that the metric
$g_2$ in the conditions of Theorem 2 is asymptotically dominated by
the metric
$$
ds^2 = g_2 = 2 d\eta \left( d\xi - \left( c_1^2 + c_2^2 \rho^2
\right) d\eta \right) + d \rho^2.
$$
Indeed, let $\gamma(s)=(\xi(s),\eta(s),p(s))$ be a causal curve in
$N$. Let $\delta_s (t)$, $0\leq t \leq \rho(p(s))$ be the shortest
normal geodesic joining the points $p_0=\delta_s(0)$ and
$p(s)=\delta_s(\rho(p(s)))$. Let
$$
\frac{d p}{ds} = d \rho \left(\frac{d p}{ds}\right)
\frac{\partial}{\partial \rho} + \left( \frac{d p}{ds} - d \rho
\left(\frac{d p}{ds}\right) \frac{\partial}{\partial \rho}\right)
$$
be the decomposition of the tangent vector into radial and tangential components
with respect to the level hypersurface of the function $\rho$ (which is, in fact, the geodesic sphere of radius $\rho$ centered at $p_0$).
Moreover, we can interpret $\frac{\partial}{\partial \rho}$ as a
tangent vector field to the normal geodesic $\delta_s$. In this case
$$
\left|\frac{d p}{ds}\right| \geq \left| d \rho \left(\frac{d
p}{ds}\right) \right| = \frac{d \rho(p(s))}{ds},\eqno{(7)}
$$
By the assumptions of the theorem, each piecewise $C^1$-smooth curve in $N$ is
projected under the map $(\xi, \eta, p) \mapsto (\xi, \eta,
\rho(p))$ to a piecewise $C^1$-smooth curve in $\R^3$. Then, (7)
implies that each causal curve in $N$ is projected under the same
map to a causal curve in $\R^3$ with the Lorentzian metric
$$
g = 2 d\eta \left( d\xi - \left( c_1^2 + c_2^2 \rho^2 \right) d\eta
\right) + d \rho^2.
$$
Obviously, the level sets of the function $\rho$ (that is, the geodesic
spheres in $M$) are compact. It remains to apply Theorem 3 and
Lemma 1.

\end{document}